\documentclass[11pt]{amsart}
\usepackage{amssymb, latexsym}
\theoremstyle{plain}
\newtheorem{theorem}{Theorem}
\newtheorem{corollary}{Corollary}

\newtheorem*{2'}{Theorem 2'}
\newtheorem*{3'}{Theorem 3'}

\theoremstyle{remark}
\newtheorem*{remark}{Remark}
\newtheorem*{Remark 1}{Remark 1}
\newtheorem*{Remark 2}{Remark 2}
\newtheorem*{Remark 3}{Remark 3}
\newtheorem*{Remark 4}{Remark 4}

\numberwithin{equation}{section}

\begin{document}

\title [Clustering in permutations avoiding patterns]
{   Clustering of consecutive numbers in permutations avoiding a pattern and in separable permutations}

\author{Ross G. Pinsky}


\address{Department of Mathematics\\
Technion---Israel Institute of Technology\\
Haifa, 32000\\ Israel}
\email{ pinsky@math.technion.ac.il}

\urladdr{http://www.math.technion.ac.il/~pinsky/}

\date{}

\begin{abstract}

Let  $S_n$ denote the set of permutations of $[n]:=\{1,\cdots, n\}$, and  denote a permutation $\sigma\in S_n$ by $\sigma=\sigma_1\sigma_2\cdots \sigma_n$.
For $l\in\{2,\cdots, n-1\}$ and $k\in\mathbb{N}$ satisfying $k+l-1\le n$, let
 $A^{(n)}_{l;k}\subset S_n$ denote the event that
 the set of  $l$ consecutive numbers $\{k, k+1,\cdots, k+l-1\}$ appears
in a  set  of consecutive positions:
$\{k,k+1,\cdots, k+l-1\}=\{\sigma_a,\sigma_{a+1},\cdots, \sigma_{a+l-1}\}$, for some $a$.
For $\tau\in S_m$, let $S_n(\tau)$ denote the set of $\tau$-avoiding permutations in
$S_n$, and let $P_n^{\text{av}(\tau)}$ denote the uniform probability measure
on $S_n(\tau)$.
Also, let $S_n^{\text{sep}}$ denote the set of separable permutations in $S_n$, and let
$P_n^{\text{sep}}$ denote the uniform probability measure on $S_n^{\text{sep}}$.
We investigate the quantities
$P_n^{\text{av}(\tau)}(A^{(n)}_{l;k})$
and  $P_n^{\text{sep}}(A^{(n)}_{l;k})$
for fixed $n$, and the limiting behavior
as $n\to\infty$.
We also consider the asymptotic properties of this limiting behavior as $l\to\infty$.

\end{abstract}

\maketitle
\section{Introduction and Statement of Results}

Let  $S_n$ denote the set of permutations of $[n]:=\{1,\cdots, n\}$, and  denote a permutation $\sigma\in S_n$ by $\sigma=\sigma_1\sigma_2\cdots \sigma_n$.
For $l\in\{2,\cdots, n-1\}$ and $k\in\mathbb{N}$ satisfying $k+l-1\le n$,  the set of  $l$ consecutive numbers $\{k,k+1,\cdots, k+l-1\}\subset [n]$ appears
in a  set  of consecutive positions in the permutation if there exists an $a$ such that $\{k,k+1,\cdots, k+l-1\}=\{\sigma_a,\sigma_{a+1},\cdots, \sigma_{a+l-1}\}$.
Let $A^{(n)}_{l;k}\subset S_n$ denote the event that
 the set of  $l$ consecutive numbers $\{k, k+1,\cdots, k+l-1\}$ appears
in a  set  of consecutive positions.
It is immediate that for any $1\le k,a\le n-l+1$,
the number of permutations $\sigma\in S_n$ satisfying  $\{k,k+1,\cdots, k+l-1\}=\{\sigma_a,\sigma_{a+1},\cdots, \sigma_{a+l-1}\}$ is equal to $l!(n-l)!$.
Thus, under the uniform probability $P_n$ on $S_n$, one has
\begin{equation}\label{uniformprob}
P_n(A^{(n)}_{l;k})=(n-l+1)\frac{l!(n-l)!}{n!}\sim\frac{l!}{n^{l-1}},\ \text{as}\ n\to\infty, \text{for}\ l\ge2.
\end{equation}
Let $A^{(n)}_l=\cup_{k=1}^{n-l+1}A^{(n)}_{l;k}$ denote the event that there exists a set of   $l$ consecutive numbers
appearing in a  set of consecutive positions.
Using inclusion-exclusion along with \eqref{uniformprob},
 it is not hard  to show that
\begin{equation*}\label{uniformtotalprob}
P_n(A^{(n)}_l)\sim\frac{l!}{n^{l-2}},\ \text{as}\ n\to\infty,\ \text{for}\ l\ge3.
\end{equation*}
In particular, the probability that a uniformly chosen permutation in $S_n$  possesses  a  cluster of three consecutive numbers in three consecutive positions
decreases to 0 as $n\to\infty$.

In this note, we study the clustering phenomenon  for permutations avoiding a fixed pattern,
as well as for separable permutations. We define and elaborate
on pattern-avoiding permutations  and separable permutations below.
Whereas $|S_n|$ increases super-exponentially, for any fixed pattern the number of permutations in $S_n$ avoiding
the pattern increases only exponentially, and so does the number of separable permutations in $S_n$. Therefore, one may expect
that the probability of $A^{(n)}_{l;k_n}$ under the uniform measure on permutations in $S_n$ avoiding a fixed
pattern and under the uniform measure on separable permutations in $S_n$ remains bounded away from zero as $n\to\infty$. For fixed $n$,
we obtain estimates on the  probability of $A^{(n)}_{l;k}$ for many patterns, and obtain exact results
for certain patterns. We also obtain exact results for separable permutations.
We  prove that the limiting  probability as $n\to\infty$ remains bounded away from zero for many patterns and for separable permutations;
whether this holds for all patterns we  leave  as
 an open question. We also focus on the asymptotic behavior of this
limiting probability as a function of $l$ as  $l\to\infty$.  With regard to the limiting probability as $n\to\infty$,
as well as the asymptotic behavior of this limit as $l\to\infty$,  the choice
of $\{k_n\}$ sometimes comes into play.

We recall the definition of pattern avoidance for permutations.
If $\sigma=\sigma_1\sigma_2\cdots\sigma_n\in S_n$ and $\tau=\tau_1\cdots\tau_m\in S_m$, where $2\le m\le n$,
then we say that $\sigma$ contains $\tau$ as a pattern if there exists a subsequence $1\le i_1<i_2<\cdots<i_m\le n$ such
that for all $1\le j,k\le m$, the inequality $\sigma_{i_j}<\sigma_{i_k}$ holds if and only if the inequality $\tau_j<\tau_k$ holds.
If $\sigma$ does not contain $\tau$, then we say that $\sigma$ \it avoids\rm\ $\tau$.
We denote by $S_n(\tau)$ the set of permutations in $S_n$ that avoid $\tau$.
If $n<m$, we define $S_n(\tau)=S_n$. We denote by
$P_n^{\text{av}(\tau)}$ the uniform probability measure on $S_n(\tau)$.

A \it separable\rm\ permutation is a permutation that can be constructed from the  singleton in $S_1$ via a series
of  iterations of
\it direct sums \rm\ and \it skew sums\rm. (See \cite{P21}, for example,  for more details.) An equivalent definition
of a separable permutation \cite{BBL} is a permutation that avoids the two patterns 2413 and 3142.
Let $S_n^{\text{sep}}$ denote the set of separable permutations in $S_n$ and let
$P_n^{\text{sep}}$ denote the uniform probability measure on $S_n^{\text{sep}}$.

We now recall several definitions and results that are required in order to state and to apply Theorem \ref{generalresult} below, which treats $P_n^{\text{av}(\tau)}(A^{(n)}_{l;k_n})$.
The permutation $\tau\in S_m$ contains the permutation  $\nu\in S_j$ \it tightly\rm\ if
for some $i$ and some $h$, one has $\tau_{i+a}=h+\nu_{1+a},\ a=0,\cdots, j-1$.  In particular,
$\tau$  contains 12 (21) tightly if $\tau_i=h$ and $\tau_{i+1}=h+1$ ($\tau_i=h+1$ and $\tau_{i+1}=h$), for some $i$ and some $h$.
One formulation of the Stanley-Wilf conjecture, completely proved in \cite{MT}, states that for every permutation $\tau\in S_m$, $m\ge2$, there exists a number
$L(\tau)$ such that
\begin{equation}\label{SW}
\lim_{n\to\infty}|S_n(\tau)|^\frac1n=L(\tau).
\end{equation}
We refer to $L(\tau)$ as the Stanley-Wilf limit.
Furthermore, it has been shown in \cite{AM} that
\begin{equation}\label{ratio}
\lim_{n\to\infty}\frac{|S_{n+1}(\tau)|}{|S_n(\tau)|}=L(\tau),
\end{equation}
if $\tau$ satisfies at least one of the following three conditions:

\begin{equation}\label{threeconditions}
\begin{aligned}
&\text{C1}.\ \tau\in S_m  \ \text{does not contain tightly at least one of 12 and 21};\\
& C2.\ \tau_1\in\{1,m\}\ \text{or}\ \tau_m\in\{1,m\};\\
& C3.\ m\ge6, \{\tau_1,\tau_2,\tau_{m-1},\tau_m\}\ \text{is a set of four consecutive numbers, and}\\
&\tau\ \text{contains exactly one tight occurrence of 12 and one tight occurrence of 21,}\\
&\text{which are in positions 1 and 2, and positions}\ m-1\ \text{and}\ m,\ \text{in either order}.
\end{aligned}
\end{equation}

We will call a permutation $\tau\in S_m$ \it cluster-free\rm\ if $\tau\not\in A^{(m)}_{l;k}$, for all $l\in\{2,\cdots, m-1\}$ and all $k\in\mathbb{N}$ satisfying
$k+l-1\le m$. (Equivalently, $\tau$ is cluster-free if and only if $\{\tau_a,\cdots,\tau_{a+l-1}\}$ is not equal to a block of $l$ consecutive numbers in $[m]$, for all
$l\in\{2,\cdots, m-1\}$ and all $a\in\mathbb{N}$ satisfying $a+l-1\le m$.)

We will prove the following theorem.
\begin{theorem}\label{generalresult}
\noindent i.  For all $\tau\in S_m$,
\begin{equation*}\label{upperone}
P_n^{\text{\rm av}(\tau)}(A^{(n)}_{l;k})\le\frac{|S_{n-l+1}(\tau)|\thinspace|S_l(\tau)|}{|S_n(\tau)|}.
\end{equation*}
ii. Assume that $\tau\in S_m$ is cluster-free. Then
  \begin{equation*}
P_n^{\text{\rm av}(\tau)}(A^{(n)}_{l;k})=\frac{|S_{n-l+1}(\tau)|\thinspace|S_l(\tau)|}{|S_n(\tau)|}.
\end{equation*}
iii. Assume that $\tau\in S_m$  does not contain tightly at least one of 12 and 21.
Then
\begin{equation*}\label{lowerone}
P_n^{\text{\rm av}(\tau)}(A^{(n)}_{l;k})\ge \frac{|S_{n-l+1}(\tau)|}{|S_n(\tau)|}.
\end{equation*}
If in fact $\tau$ contains tightly neither 12 nor 21, then
\begin{equation*}\label{lowerone}
P_n^{\text{\rm av}(\tau)}(A^{(n)}_{l;k})\ge \frac{2|S_{n-l+1}(\tau)|}{|S_n(\tau)|}.
\end{equation*}
\end{theorem}
If a permutation is cluster-free, then it satisfies  C1 in \eqref{threeconditions}.
Theorem \ref{generalresult} along with \eqref{ratio} and \eqref{threeconditions} immediately yield the following corollary.
\begin{corollary}\label{Cor1}
Let $L(\tau)$ denote the Stanley-Wilf limit.
Let $\{k_n\}$ be an arbitrary sequence satisfying $1\le k_n\le n-l+1$.

\noindent i.  Assume that $\tau\in S_m$ satisfies at least one of the conditions in \eqref{threeconditions}.
 Then
\begin{equation}\label{upperbound}
\limsup_{n\to\infty}P_n^{\text{\rm av}(\tau)}(A^{(n)}_{l;k_n})\le\frac{|S_l(\tau)|}{(L(\tau))^{l-1}}.
\end{equation}
ii. Assume that $\tau\in S_m$ is cluster-free. Then
\begin{equation}\label{exactlimit}
\lim_{n\to\infty}P_n^{\text{\rm av}(\tau)}(A^{(n)}_{l;k_n})=\frac{|S_l(\tau)|}{(L(\tau))^{l-1}}.
\end{equation}
In particular
\begin{equation*}
\lim_{n\to\infty}P_n^{\text{\rm av}(\tau)}(A^{(n)}_{l;k_n})=\begin{cases}\frac{l!}{(L(\tau))^{l-1}},\ \text{if}\ l<m;\\
\frac{m!-1}{(L(\tau))^{m-1}},\ \text{if}\ l=m.\end{cases}
\end{equation*}
iii.  Assume that $\tau\in S_m$  does not contain tightly at least one of 12 and 21. Then
\begin{equation}\label{lowerbound}
\liminf_{n\to\infty}P_n^{\text{\rm av}(\tau)}(A^{(n)}_{l;k_n})\ge\frac1{(L(\tau))^{l-1}}.
\end{equation}
If in fact $\tau$ contains tightly neither 12 nor 21, then
\begin{equation*}
\liminf_{n\to\infty}P_n^{\text{\rm av}(\tau)}(A^{(n)}_{l;k_n})\ge\frac2{(L(\tau))^{l-1}}.
\end{equation*}
\end{corollary}
Each of the three remarks below deals with a different part of Corollary \ref{Cor1}.
\begin{remark}  It follows from
\eqref{lowerbound} that
$\liminf_{n\to\infty}P_n^{\text{av}(\tau)}(A^{(n)}_{l;k_n})>0$, if $\tau$ does not contain tightly at least one of 12 and 21.
We pose the following question:
\end{remark}

\bf\noindent Question.\rm\ Is
$\liminf_{n\to\infty}P_n^{\text{av}(\tau)}(A^{(n)}_{l;k_n})>0$, for all $\tau\in \cup_{m=2}^\infty S_m$,
for all $\{k_n\}$ and for all $l\ge2$?
\medskip

\begin{remark}
To the author's knowledge, in all of the small number of cases where the asymptotic behavior of
$|S_l(\tau)|$ as $l\to\infty$  is known, one has $\lim_{l\to\infty}\frac{|S_l(\tau)|}{(L(\tau))^l}=0$.
This is true for all $\tau\in S_3$---see below.
For $\tau=12\cdots m$,  one has $|S_n(12\cdots m)|\sim c(m-1)^{2n}n^{-\frac12(m^2-2m)}$,
for a constant $c$ given by a multiple integral \cite{R}.
For $\tau=1342$, one has
$|S_n(1342)|\sim \frac6{243\sqrt\pi}8^nn^{-\frac52}$ \cite{B97}.
From  \eqref{upperbound},
$\lim_{l\to\infty}\frac{|S_l(\tau)|}{(L(\tau))^l}=0$
is  a sufficient condition for
$\lim_{l\to\infty}\limsup_{n\to\infty}P_n^{\text{av}(\tau)}(A^{(n)}_{l;k_n})=0$, if
$\tau$ satisfies one of the conditions in \eqref{threeconditions}.
When a sequence of  probability measures $\{\mathcal{P}_n\}$  on the sequence of spaces $\{S_n\}$ is such that
 $\lim_{l\to\infty}\limsup_{n\to\infty} \mathcal{P}_n(A^{(n)}_{l;k_n})>0$, for some choice of $\{k_n\}_{n=1}^\infty$, we say that super-clustering occurs; see \cite{P} for situations where this
phenomenon is encountered. We make the following conjecture.
\end{remark}

\noindent \bf Conjecture.\rm\
$\lim_{l\to\infty}\limsup_{n\to\infty}P_n^{\text{av}(\tau)}(A^{(n)}_{l;k_n})=0$,
for all $\tau\in \cup_{m=2}^\infty S_m$ and
for all $\{k_n\}$.

\begin{remark}
From \eqref{exactlimit} and \eqref{SW}, it follows that if $\tau$ is cluster-free, then $\lim_{n\to\infty}P_n^{\text{\rm av}(\tau)}(A^{(n)}_{l;k_n})$ exhibits sub-exponential
decay in the length $l$ of the cluster, for all choices of $\{k_n\}_{n=1}^\infty$. For a situation with exponential decay, see Remark \ref{Rem123}.
\end{remark}

\medskip

We now turn to the separable permutations $S^{\text{sep}}_n$.
\begin{theorem}\label{separablethm}
 \begin{equation*}
P_n^{\text{\rm sep}}(A^{(n)}_{l;k})=\frac{|S^{\text{\rm sep}}_{n-l+1}|\thinspace|S_l^{\text{\rm sep}}|}{|S_n^\text{\rm sep}|}.
\end{equation*}
\end{theorem}
\medskip

The generating function for the enumeration of separable permutations is known explicitly and allows one
to show \cite[p. 474-475]{FS} that
\begin{equation}\label{sepasymp}
|S_n^{\text{\rm sep}}|\sim\frac1{2\sqrt{\pi n^3}}(3-2\sqrt2)^{-n+\frac12}.
\end{equation}
Theorem \ref{separablethm} and \eqref{sepasymp} immediately yield the following corollary.
\begin{corollary}\label{sepcor}
Let $\{k_n\}$ be an arbitrary sequence satisfying $1\le k_n\le n-l+1$.
Then
\begin{equation*}
\lim_{n\to\infty}P_n^{\text{\rm sep}}(A^{(n)}_{l;k_n})=(3-2\sqrt2)^{l-1}|S_l^{\text{\rm sep}}|.
\end{equation*}
In particular,
$$
\begin{aligned}
&\lim_{n\to\infty}P_n^{\text{\rm sep}}(A^{(n)}_{3;k_n})=6(3-2\sqrt2)^2;\\
&\lim_{n\to\infty}P_n^{\text{\rm sep}}(A^{(n)}_{4;k_n})=22(3-2\sqrt2)^3.
\end{aligned}
$$
\end{corollary}
\begin{remark} Note that as the length $l$ of the cluster grows, the
quantity $\lim_{n\to\infty}P_n^{\text{\rm sep}}(A^{(n)}_{l;k_n})$
decays on the order $l^{-\frac32}$.
\end{remark}

We now return to permutations avoiding one pattern and  consider in additional detail the case  $m=3$.
 Note that all $\tau$ in $S_3$ satisfy the condition
in part (iii) of Theorem \ref{generalresult} and the conditions in parts (i) and (iii) of Corollary  \ref{Cor1};
thus the theorem and  corollary give upper and lower bounds
on $P_n^{\text{av}(\tau)}(A^{(n)}_{l;k})$ and
$\lim_{n\to\infty}P_n^{\text{av}(\tau)}(A^{(n)}_{l;k_n})$, for all $\tau$ in $S_3$.
However,  no $\tau$ in $S_3$ satisfies the condition in part (ii) of the theorem and corollary.
The result below for the case $m=3$  reveals some  phenomena that cannot be gleaned from
Theorem \ref{generalresult} and Corollary  \ref{Cor1}.

We recall the well known fact \cite{B} that $|S_n(\tau)|=C_n$, for all
$\tau\in S_3$, where $C_n$ is
the $n$th Catalan number, which satisfies
\begin{equation}\label{catalan}
C_n=\frac1{n+1}\binom{2n}n,\ n\ge0, \text{and} \ C_n\sim  \frac{4^n}{\sqrt\pi \thinspace n^\frac32}\ \text{as}\ n\to\infty.
\end{equation}

\begin{theorem}\label{m=3}
 Let $\tau\in\{123,321\}$. Then
 \begin{equation}\label{exact321}
 P_n^{\text{\rm av}(\tau)}(A^{(n)}_{l;k})=\frac{C_{n-l+1}+C_{k-1}C_{n-k-l+1}(C_l-1)}{C_n}.
 \end{equation}
\end{theorem}
From \eqref{catalan} and \eqref{exact321}, a straight forward calculation yields the following corollary.
\begin{corollary}\label{Cor2}
Let $\tau\in\{123,321\}$.

\noindent i.
\begin{equation*}
\lim_{n\to\infty}P_n^{\text{\rm av}(\tau)}(A^{(n)}_{l;k})=
\lim_{n\to\infty}P_n^{\text{\rm av}(\tau)}(A^{(n)}_{l;n+2-k-l})=\frac1{4^{l-1}}+\frac{C_{k-1}(C_l-1)}{4^{k+l-1}}.
\end{equation*}
ii. If $\lim_{n\to\infty}k_n=\lim_{n\to\infty}(n-k_n)=\infty$, then
\begin{equation*}
\lim_{n\to\infty}P_n^{\text{\rm av}(\tau)}(A^{(n)}_{l;k_n})=\frac1{4^{l-1}}.
\end{equation*}
\end{corollary}
\begin{remark}\label{Rem123} Note that for $\tau\in\{123,321\}$ and  $\lim_{n\to\infty}k_n=\lim_{n\to\infty}(n-k_n)=\infty$,  the quantity $\lim_{n\to\infty}P_n^{\text{av}(\tau)}(A^{(n)}_{l;k_n})$
decays exponentially in the length $l$ of the cluster. On the other hand,
$\lim_{n\to\infty}P_n^{\text{av}(\tau)}(A^{(n)}_{l;k})$ and  $\lim_{n\to\infty}P_n^{\text{av}(\tau)}(A^{(n)}_{l;n+2-k-l})$
decay  on the order $l^{-\frac32}$ as $l\to\infty$, and on the order $(lk)^{-\frac32}$ as $k,l\to\infty$.

\end{remark}

In section \ref{GRproof} we prove Theorems \ref{generalresult} and \ref{separablethm}, and  in section \ref{thmproof} we prove Theorem
\ref{m=3}.

\section{Proofs of Theorems \ref{generalresult} and \ref{separablethm}}\label{GRproof}

\noindent\it Proof of Theorem \ref{generalresult}.\rm\
Fix $\tau$, $k$ and $l$ as in the statement of the theorem.
For $a\in\{1,\cdots, n-l+1\}$, define
$$
A^{(n)}_{l;k;a}=\big\{\sigma\in A^{(n)}_{l;k}:\{k,k+1,\cdots, k+l-1\}=\{\sigma_a,\sigma_{a+1},\cdots, \sigma_{a+l-1}\}\big\}.
$$
Then the sets $\{A^{(n)}_{l;k;a}\}_{a=1}^{n-l+1}$ are disjoint and
$A^{(n)}_{l;k}=\cup_{a=1}^{n-l+1}A^{(n)}_{l;k;a}$.

If $\nu=\{\nu_i\}_{i=1}^{|B|}$ is a permutation of a finite set  $B\subset\mathbb{N}$, let
$\nu^{B^{-1}}$ denote the permutation it naturally induces on $S_{|B|}$.
Thus, for example, if $B=\{3,5,7,8\}$ and $\nu=7385$, then $\nu^{B^{-1}}=3142\in S_4$.
Conversely, if  $\nu$ is a permutation of $S_{|B|}$, let   $\nu^{B}$ denote
the permutation it naturally induces on $B$.

Until further notice,  consider $a$ fixed.
Let $\sigma\in A^{(n)}_{l;k;a}\cap S_n(\tau)$.
We  describe a procedure    to contract $\sigma$ to a permutation in $S_{n-l+1}(\tau)$.
Define the permutation $\overline\sigma=\overline\sigma(\sigma)$ of the set
$B=\{1,\cdots, k,k+l,\cdots, n\}$  by
$$
\overline\sigma_i=\begin{cases}\sigma_i,\ 1\le i\le a-1;\\ k,\ i=a;\\ \sigma_{i+l-1}, \ i=a+1,\cdots, n-l+1.\end{cases},
$$
and define
$$
\eta=\eta(\sigma)=\overline\sigma^{B^{-1}}  \in S_{n-l+1}.
$$
It  follows from the construction  that
\begin{equation}\label{eta}
\eta\in S_{n-l+1}(\tau) \ \text{and}\ \eta_a=k.
\end{equation}

We concretize the above construction with an example.  Let $n=9, l=3, k=a=4$ and $\tau=123\in S_3$.
Let $\sigma=798645312\in A^{(9)}_{3;4;4}\cap S_9(123)$. The set $B$ is given by $B=\{1,2,3,4,7,8,9\}$ and $\bar\sigma=7984312$---the cluster $645$  in $\sigma$ has been
contracted to $4$ in $\bar\sigma$. Finally, $\eta=\eta(\sigma)=\bar\sigma^{B^{-1}}=5764312$ satisfies $\eta\in S_7(123)$ and $\eta_4=4$.

Obviously the map taking $\sigma\in A^{(n)}_{l;k;a}\cap S_n(\tau)$ to $\eta(\sigma)$ is not injective.
However,
\begin{equation}\label{partialinject}
\begin{aligned}
&\eta(\sigma)\neq\eta(\sigma'), \ \text{if}\ \sigma,\sigma'\in A^{(n)}_{l;k;a}\cap S_n(\tau)\ \text{are distinct and satisfy}\\
&\sigma_{a+i}=\sigma'_{a+i}, \ i=0,\cdots, l-1.
\end{aligned}
\end{equation}

Conversely, let $\eta$ satisfy \eqref{eta}. We describe a procedure to extend
$\eta$  to a permutation in $A^{(n)}_{l;k;a}$, which may or may not belong to $S_n(\tau)$.
Let $B=\{1,\cdots, k,k+l,\cdots, n\}$ as above.
For each $\rho\in S_l$, define
 $\sigma^{\rho}=\sigma^{\rho}(\eta)\in S_n$ by
\begin{equation}\label{etasigma}
\sigma_i^{\rho}=\begin{cases}\eta^B_i,\ i=1,\cdots, a-1;\\ k-1+\rho_{i-a+1},\ i=a,\cdots, a+l-1;\\
\eta^B_{i-l+1},\ i=a+l,\cdots, n.
\end{cases}
\end{equation}
It follows from the construction that
\begin{equation}\label{sigmarho}
\sigma^\rho(\eta)\in A^{(n)}_{l;k;a}.
\end{equation}
Also, of course, the map taking $\eta$ satisfying \eqref{eta} to $\sigma^\rho(\eta)$ is injective.

As an example of the above construction, again with $n=9, l=3, k=a=4$ and $\tau=123\in S_3$, let $\eta=5764312$. Then $\eta$ satisfies  \eqref{eta}.
We have $B=\{1,2,3,4,7,8,9\}$. Choose, for example, $\rho=213\in S_3$. Then $\eta^B=7984312$ and $\sigma^\rho=\sigma^\rho(\eta)=798546312\in S_9$---the $4$ in
$\eta^B$ has been expanded to the cluster $546$ in $\sigma^\rho$.

We now consider the question of whether or not $\sigma^\rho\in S_n(\tau)$.
Since $\eta$ is $\tau$-avoiding, $\eta^B$ is also $\tau$-avoiding. The transition
from $\eta^B$ to $\sigma^\rho$ in \eqref{sigmarho} expands the $a$th position in $\eta^B$, which contains the number
$k$, into $l$ positions which contain the numbers $\{k,\cdots,k+l-1\}$ in the order dictated by the pattern $\rho$.
Of course, a necessary condition for $\sigma^\rho\in S_n(\tau)$ is that
$\rho\in S_l(\tau)$.

Assume now that
$\tau$ does not contain tightly 12. If  the numbers
$\{k,\cdots, k+l-1\}$  are placed in ascending  order, that is, if we choose $\rho=1\cdots l$, then the resulting permutation $\sigma^\rho$ retains
the $\tau$-avoiding property.
Indeed, assume to the contrary that  the pattern $\tau$  appears in $\sigma^\rho$. Then this $\tau$ pattern will include at least two of the numbers $\{k,\cdots, k+l-1\}$ that appear in increasing order in the positions
$\{a,\cdots, a+l-1\}$ of $\sigma^\rho$ (because otherwise the pattern $\tau$ would already   appear in $\eta^B$). But then $\tau$ will contain tightly 12, which is a contradiction.
Similarly, if we assume that  $\tau$ does not contain tightly 21, and the numbers
$\{k,\cdots, k+l-1\}$  are placed in descending  order, that is, we choose $\rho=l\cdots 1$, then the resulting permutation $\sigma^\rho$ retains
the $\tau$-avoiding property.

Now assume that  $\tau$ is cluster-free.  Let  the numbers $\{k,\cdots, k+l-1\}$ be placed in any order
that is $\tau$-avoiding, that is  choose any $\rho\in S_l(\tau)$. Then the
resulting permutation $\sigma^\rho$ retains
the $\tau$-avoiding property.
Indeed, assume to the contrary that  the pattern $\tau$  appears in $\sigma^\rho$. Then this $\tau$ pattern will include at least two of the numbers $\{k,\cdots, k+l-1\}$ that appear in the positions
$\{a,\cdots, a+l-1\}$ of $\sigma^\rho$ (because otherwise the pattern $\tau$ would  already appear in $\eta^B$).
 This $\tau$ pattern will also include at least  one of the numbers in $[n]-\{k,\cdots, k+l-1\}$
that appear in the positions  $[n]-\{a,\cdots, a+l-1\}$ of $\sigma^\rho$ (because otherwise the pattern  $\tau$ would appear in $\rho$).
For definiteness, let $r$ denote  the number of numbers from  $\{k,\cdots, k+l-1\}$ that appear in this $\tau$ pattern.
Then $2\le r\le \min(m-1,l)$.
It  then  follows that $\tau$ has a cluster of size $r$, which contradicts the assumption that $\tau$ is cluster-free.

With the above constructions and commentary, we can  prove the theorem.
We begin with part (iii). Assume that $\tau\in S_m$ does not contain tightly at   least one of
12 and 21. If $\tau$ contains tightly 12, let $\rho=l\cdots1\in S_l$, if $\tau$ contains tightly 21, let
$\rho=1\cdots l\in S_l$, and if $\tau$ contains tightly neither 12 nor 21, let $\rho$ be
either $1\cdots l$ or $l\cdots 1$.
Then for each $\eta$ satisfying \eqref{eta}, it follows from the previous paragraph that
$\sigma^\rho(\eta)\in S_n(\tau)$.
This in conjunction with  \eqref{sigmarho} gives
\begin{equation}\label{keypartiii}
\sigma^\rho(\eta)\in A^{(n)}_{l;k;a}\cap S_n(\tau),\ \text{for}\ \eta\ \text{satisfying}\ \eqref{eta}.
\end{equation}
From \eqref{keypartiii} and the fact that the map $\eta\to\sigma^\rho(\eta)$ is injective, it follows
that
\begin{equation}\label{setinequalalower}
|A^{(n)}_{l;k;a}\cap S_n(\tau)|\ge|\{\eta\in S_{n-l+1}(\tau):\eta_a=k\}|.
\end{equation}
Considering  the union over $a\in\{1,\cdots, n-l+1\}$ of the sets on the left and right hand sides
of \eqref{setinequalalower}, we obtain
\begin{equation}\label{union}
|A^{(n)}_{l;k}\cap S_n(\tau)|\ge|S_{n-l+1}(\tau)|.
\end{equation}
From \eqref{union} we conclude that
\begin{equation*}
P_n^{\text{av}(\tau)}(A^{(n)}_{l;k})=\frac{|A^{(n)}_{l;k}\cap S_n(\tau)|}{|S_n(\tau)|}\ge\frac{|S_{n-l+1}(\tau)|}{|S_n(\tau)|}.
\end{equation*}
This proves part (iii) in the case that $\tau$ does not contain tightly at least one of 12 and 21.
If in fact $\tau$ contains tightly neither 12 nor 21, then we can use both $\rho=1\cdots l$ and $\rho=l\cdots 1$ in the construction above, which then gives
\begin{equation*}
P_n^{\text{av}(\tau)}(A^{(n)}_{l;k})=\frac{|A^{(n)}_{l;k}\cap S_n(\tau)|}{|S_n(\tau)|}\ge\frac{2|S_{n-l+1}(\tau)|}{|S_n(\tau)|}.
\end{equation*}
This completes the  proof of part (iii).

We now turn to the proof of part (i). Let $\tau\in S_m$ be arbitrary.
Let $\sigma\in A^{(n)}_{l;k;a}\cap S_n(\tau)$. For any $\rho\in S_l$, map $\sigma$
to $\eta=\eta(\sigma)$,
and then map  $\eta(\sigma)$ to $\sigma^\rho(\eta(\sigma))$.
As noted above, a necessary condition for  $\sigma^\rho(\eta(\sigma))$ to avoid $\tau$ is that
$\rho\in S_l(\tau)$.
Since the map $\eta$ satisfies
 \eqref{partialinject}, we conclude that
 \begin{equation}\label{setinequalaupper}
|A^{(n)}_{l;k;a}\cap S_n(\tau)|\le|\{\eta\in S_{n-l+1}(\tau):\eta_a=k\}|\thinspace|S_l(\tau)|.
\end{equation}
Considering  the union over $a\in\{1,\cdots, n-l+1\}$ of the sets on the left and right hand sides
of \eqref{setinequalaupper}, we obtain
\begin{equation}\label{union2}
|A^{(n)}_{l;k}\cap S_n(\tau)|\le|S_{n-l+1}(\tau)|\thinspace|S_l(\tau)|.
\end{equation}
From \eqref{union2} we conclude that
\begin{equation*}
P_n^{\text{av}(\tau)}(A^{(n)}_{l;k})=\frac{|A^{(n)}_{l;k}\cap S_n(\tau)|}{|S_n(\tau)|}\le\frac{|S_{n-l+1}(\tau)|\thinspace|S_l(\tau)|}{|S_n(\tau)|},
\end{equation*}
which proves part (i).

We now prove part (ii).
Assume that  $\tau\in S_m$ is cluster-free.
Then as noted above, $\sigma^\rho\in S_n(\tau)$, for all $\rho\in S_l(\tau)$, and
$\sigma^\rho\not\in S_n(\tau)$
for all $\rho\not\in S_l(\tau)$.
Thus, from the considerations invoked in the proofs of parts (iii) and (i),
it follows that
 \begin{equation}\label{setinequalaequal}
|A^{(n)}_{l;k;a}\cap S_n(\tau)|=|\{\eta\in S_{n-l+1}(\tau):\eta_a=k\}|\thinspace|S_l(\tau)|.
\end{equation}
Considering  the union over $a\in\{1,\cdots, n-l+1\}$ of the sets on the left and right hand sides
of \eqref{setinequalaequal}, we obtain
\begin{equation}\label{union3}
|A^{(n)}_{l;k}\cap S_n(\tau)|=|S_{n-l+1}(\tau)|\thinspace|S_l(\tau)|.
\end{equation}
From \eqref{union3} we conclude that
\begin{equation*}
P_n^{\text{av}(\tau)}(A^{(n)}_{l;k})=\frac{|A^{(n)}_{l;k}\cap S_n(\tau)|}{|S_n(\tau)|}=\frac{|S_{n-l+1}(\tau)|\thinspace|S_l(\tau)|}{|S_n(\tau)|},
\end{equation*}
which proves part (ii).
\hfill $\square$

\medskip

\it\noindent Proof of Theorem \ref{separablethm}.\rm\ For permutations $\tau, \eta\in S_m$, let $S_n(\tau,\eta)$ denote the set of permutations in $S_n$ that avoid
both $\tau$ and $\eta$, and let $P_n^{\text{av}(\tau,\eta)}$ denote the uniform probability measure on $S_n(\tau,\eta)$.
It is easy to check that the following analog  of part (ii) of Theorem \ref{generalresult} can be proved in exactly the same way as its antecedent was proven:
\begin{equation}\label{iiclusterfree}
P_n^{\text{av}(\tau,\eta)}(A^{(n)}_{l;k}))=\frac{|S_{n-l+1}(\tau,\eta)|\thinspace|S_l(\tau,\eta)|}{|S_n(\tau,\eta)|},\ \text{if}\ \tau\ \text{and}\ \eta\ \text{are
cluster-free}.
\end{equation}
Recalling that the class of separable permutations is the class of  permutations   that avoid the patterns 2413 and 3142, and noting
that the permutations 2413 and 3142 are cluster-free, the theorem follows from \eqref{iiclusterfree}.
 \hfill $\square$

\section{Proofs of Theorem \ref{m=3}}\label{thmproof}
We will assume that $\tau=321$.  The proof for $\tau=123$ follows by making obvious modifications.
For the time being, fix $a\in\{1,\cdots, n-l+1\}$.
Recall the constructions in the proof of Theorem \ref{generalresult}, with $\tau=321$.
A permutation $\sigma\in A^{(n)}_{l;k;a}\cap S_n(321)$
is mapped to $\eta=\eta(\sigma)$ satisfying \eqref{eta}, that is
\begin{equation}\label{etaagain}
\eta\in S_{n-l+1}(321) \ \text{and}\ \eta_a=k.
\end{equation}
And for each $\rho\in S_l$, a permutation $\eta$ satisfying \eqref{etaagain} is mapped
to $\sigma^\rho=\sigma^\rho(\eta)$ satisfying
\eqref{sigmarho}, that is,
\begin{equation}\label{sigmarhoagain}
\sigma^\rho(\eta)\in A^{(n)}_{l;k;a}.
\end{equation}

We now investigate when in fact $\sigma^\rho(\eta)\in S_n(321)\cap A^{(n)}_{l;k;a}$.
If $\rho\not\in S_l(321)$, then of course $\sigma^\rho\not\in S_n(321)$.
If $\rho=12\cdots l$, then $\sigma^\rho\in S_n(321)$, for all
$\eta$ satisfying \eqref{etaagain}.
This contributes $|\{\eta\in S_{n-l+1}(321):\eta_a=k\}|$ members
to $S_n(321)\cap A^{(n)}_{l;k;a}$.

Now consider any of the other $\rho\in S_l$. There are $C_l-1$ such $\rho$.
Since $\rho$ has a decreasing subsequence of length 2, in order
to have $\sigma^\rho(\eta)\in S_n(321)$, all of the numbers $\{1,\cdots, k-1\}$ must appear among the first $a-1$ positions of $\eta$, and all the numbers $\{k+1,\cdots n-l+1\}$ must appear
 among the last  $n-a-l+1$ positions of $\eta$. This is possible only if $k=a$.
 If indeed $k=a$, then $\sigma^\rho(\eta)\in S_n(321)$ if
 and only if  the first $k-1$ positions of $\eta$ are filled in a 321-avoiding way by the numbers
$\{1,\cdots, k-1\}$  and the last  $n-k-l+1$ positions of $\eta$ are filled in a 321-avoiding
way by the numbers $\{k+1,\cdots n-l+1\}$. (The one remaining position, position $a$, is by assumption filled by the number $k$.)
Thus, the  number of such $\eta$ is $C_{k-1}C_{n-k-l+1}$.
Therefore, if $a=k$, we obtain an additional contribution of
$C_{k-1}C_{n-k-l+1}(C_l-1)$ members to $S_n(321)\cap A^{(n)}_{l;k;a}$.
We conclude from this and from the construction in the proof of Theorem \ref{generalresult} that
\begin{equation}\label{ak}
|A^{(n)}_{l;k;a}\cap S_n(321|=\begin{cases}|\{\eta\in S_{n-l+1}(321):\eta_a=k\}|,\ \text{if}\ a\neq k;\\
|\{\eta\in S_{n-l+1}(321):\eta_a=k\}|+C_{k-1}C_{n-k-l+1}(C_l-1), \ \text{if}\ a=k.\end{cases}
\end{equation}
Considering the union over all $a\in\{1,\cdots, n-l+1\}$ of the sets on the left and right hand sides
of \eqref{ak}, we obtain
\begin{equation}\label{union4}
|A^{(n)}_{l;k}\cap S_n(321)|=C_{n-l+1}+C_{k-1}C_{n-k-l+1}(C_l-1).
\end{equation}
From \eqref{union4} we conclude that
$$
P_n^{\text{av}(321)}(A^{(n)}_{l;k})=\frac{|A^{(n)}_{l;k}\cap S_n(321)|}{|S_n(321)|}=
\frac{C_{n-l+1}+C_{k-1}C_{n-k-l+1}(C_l-1)}{C_n},
$$
which proves the theorem. \hfill $\square$

\bf\noindent Acknowledgment.\rm\ The author thanks Toufik Mansour for supplying him with   the reference \cite{AM}.


\begin{thebibliography}{99}





\bibitem{AM}
 Atapour, M. and  Madras, N., \emph{Large deviations and ratio limit theorems for pattern-avoiding permutations},  Combin. Probab. Comput. \textbf{23} (2014), 161-200.
 \MR{3166466}

\bibitem{B97} Bona, M., \emph{Exact enumeration of 1342-avoiding permutations: a close link with labeled trees and planar maps}, J. Combin. Theory Ser. A \textbf{80} (1997), 257-272.
\MR{1485138}
\bibitem{B} B\'ona, M., \emph{On three different notions of monotone subsequences}, Permutation patterns, 89-114, London Math. Soc. Lecture Note Ser., \textbf{376}, Cambridge Univ. Press, Cambridge, (2010).
\MR{2732825}
\bibitem{BBL}
 Bose, P.,  Buss, J. and Lubiw, A., \emph{Pattern matching for permutations}, Inform. Process. Lett. \textbf{65} (1998), 277-283.
 \MR{1620935}


\bibitem{FS} Flajolet, P. and Sedgewick, R.,
 \emph{Analytic Combinatorics},  Cambridge University Press, Cambridge, (2009).
 \MR{2483235}


\bibitem{MT}
 Marcus, A. and  Tardos, G., \emph{Excluded permutation matrices and the Stanley-Wilf conjecture}, J. Combin. Theory Ser. A \textbf{107} (2004), 153-160.
 \MR{2063960}




\bibitem{P21} Pinsky, R. \emph{The infinite limit of separable permutations},  Random Structures Algorithms \textbf{59} (2021), 622-639,
\MR{4323312}


\bibitem{P} Pinsky, R., \emph{Clustering of consecutive numbers in permutations under Mallows distributions  and super-clustering under general $p$-shifted distributions}, to
appear in the Electronic Journal of Probability.

\bibitem{R}
 Regev, A.,  \emph{Asymptotic values for degrees associated with strips of Young diagrams}, Adv. in Math. \textbf{41} (1981),  115-136.
 \MR{0625890}

\bibitem{W} Wolfowitz, J. \emph{Note on runs of consecutive elements}, Ann. Math. Statistics \textbf{15} (1944), 97-98.
\MR{0010341}



\end{thebibliography}
\end{document}